\documentclass{amsart}
\usepackage{amscd} 
\usepackage{amsfonts} 
\usepackage{amssymb} 
\usepackage{latexsym} 
\input{diagrams}
\newcommand{\ncm}{\newcommand}

\def\C{\mathbb{C}\,}

\newtheorem{theorem}{Theorem}[section]
\newtheorem{prop}[theorem]{Proposition}
\newtheorem{lemma}[theorem]{Lemma}
\newtheorem{cor}[theorem]{Corollary}

\newtheorem{lem&def}[theorem]{Lemma \& Definition}
\newtheorem{definition}[theorem]{Definition}
\newtheorem{example}[theorem]{Example}

\ncm{\End}{\mbox{\rm End}\,}

\def\Hom{\mbox{\rm Hom}\,}

\def\Im{\mbox{\rm Im}\,}

\def\id{\mbox{\rm id}}

\def\into{\hookrightarrow}
\def\to{\rightarrow}

\def\o{\otimes}    %tensor product 
\def\b{\, \Box}    %cotensor product
\def\bra{\langle}
\def\ket{\rangle}

\ncm{\rarr}[1]{\stackrel{#1}{\longrightarrow}}
\ncm{\larr}[1]{\stackrel{#1}{\longleftarrow}}
\def\cop{\Delta}

\def\eps{\varepsilon}

\def\du1{\hat 1}

\def\-1{_{(-1)}}
\def\0{_{(0)}}
\def\1{_{(1)}}
\def\2{_{(2)}}
\def\3{_{(3)}}
\def\4{_{(4)}}

\def\|{\, | \,}

\def\du1{\hat 1}

\def\ract{\triangleleft}

\begin{document}

\title[Codepth Two and Related Topics]{Codepth Two and Related Topics}
\author{Lars Kadison}
\address{Department of Mathematics \\ University of Pennsylvania \\ 209 South 33rd Street \\ Philadelphia, PA 19104-6395} 
\email{lkadison@c2i.net} 
\date{}
\thanks{}
 \dedicatory{Dedicated to Daniel Kastler on his eightieth  birthday}
 
\begin{abstract} 
A depth two extension $A \| B$ is shown to be weak depth
two over its double centralizer $V_A(V_A(B))$ if this is separable over $B$.
We consider various examples and non-examples of depth one and two properties.
Depth two and its relationship to direct and tensor product
of algebras as well as cup product of relative Hochschild cochains is
examined.   
Section~6 introduces a notion of codepth two coalgebra homomorphism $g: C \to D$, dual to a depth two algebra homomorphism. It is shown that the endomorphism ring of bicomodule endomorphisms $\End {}^DC^D$ forms a right bialgebroid over 
the centralizer subalgebra $g^*: D^* \to C^*$ of the dual algebra $C^*$.       
\end{abstract} 
\maketitle

\section{Introduction}

From the quantum algebraic viewpoint,
a depth two subalgebra is a notion that generalizes  finite index normal Hopf subalgebra.  Somewhat like a normal Hopf subalgebra is a Hopf-Galois extension,  so does a depth two subalgebra admit a Galois theory
of bialgebroid-valued coactions.   The notion of depth two is also closely related to older definitions of normal subalgebra via representation theory, such as Rieffel's,
where normality is defined by invariance of
contracted maximal ideals with respect to the over-algebra. For example, 
a depth two extension has normal centralizer. Among Hopf subalgebras, the notion of normal subalgebra is the same as the notion of normal Hopf subalgebra.  It follows that
a depth two Hopf subalgebra has normal double centralizer subalgebra under certain circumstances
\cite{LK2005}. In Section~2 of this paper we pursue a more general fact underlying this;
 that a depth two extension $A \| B$ where the double centralizer $W = V_A(V_A(B))$
is a separable extension of $B$ has \textit{weak} depth two extension $A \| W$. The definition of weak depth two, given in section~2, 
is modelled on Mewborn and McMahon's weakening of H-separability to ``strong separability'' in \cite{MM}.  
In section~3 we consider counterexamples and
propositions for other closure and
transitivity properties of depth two in a tower of algebras $A \supset C \supset B$.
The aim is ultimately to extend Galois correspondences for special H-separable extensions
to certain depth two algebra extensions.  In section~4, we characterize a f.g.\ projective left Galois extensions with bialgebroid  action in terms of its smash product being isomorphic
to its right endomorphism ring.  In section~5, we note that a one-sided depth two extension
$A \| B$ has  relative Hochschild cochains with values in $A$ generated w.r.t.\ the cup
product by the $1$-cochains.

The notion of bialgebroid dualizes to the notion of
bicoalgebroid in \cite{BM}.  It turns out that a weak
Hopf algebra is both a bialgebroid and bicoalgebroid
over its target subalgebra.  To this author, it is an interesting continuation of this inquiry to dualize depth two algebra homomorphism, define a codepth two
coalgebra homomorphism and look for bicoalgebroids, or at least bialgebroids.  One reason one might
want to embark on this project is that a Hopf algebra homomorphism $K \to H$
is both an algebra and coalgebra homomorphism, so 
considering it as a codepth two coalgebra homomorphism is in principle just as interesting as considering it as
a depth two algebra homomorphism; the latter having
been done in terms of ``depth two Hopf subalgebras,'' the former then leading perhaps to interesting
cases of ``codepth two Hopf quotient algebras.''
We carry out the beginnings of this project in section~6, where we define a left codepth two coalgebra homomorphism in terms of its associated cotensor product.  The concept of codepth two is expressed in terms of coordinates called quasibases. It is then shown
that the endomorphism algebra of the canonical bicomodule is a bialgebroid over the centralizer of
the dual algebra homomorphism.  We conclude with a discussion
of how one might continue this investigation.

%%%%%%%%%%%%%%%%%%%%%%%%%%%%%%%%%%%%%%%%%%%%%%%%%%%%%%%%%%%%%%%%%%%%%%%%%%%%%%%%%%%%%%%%%%%%%%%%%%%%%%%%%%%%%%%%%%%%%
\section{Preliminaries and weak depth two}

In this paper, algebras are unital associative over a commutative ground ring $K$ and are themselves not necessarily commutative.
An algebra extension $A \| B$ is a unit-preserving algebra homomorphism $B \to A$,
a proper extension if this mapping is monic.  We focus mostly on the induced
bimodule ${}_BA_B$ and mostly suppress the homomorphism.  Unadorned tensors, hom-groups and endomorphism-groups between algebras are
over the ground ring unless otherwise stated. In this section we denote $A^B = \{ a \in A | \forall\, b \in B, \, ba = ab   \} $, the centralizer $V_A(B)$ of $B$ in $A$, although $A^B$ should not be confused
with the comodule notation in section~6, nor the invariant subring notation $A^S$ where
a bialgebroid $S$ acts on $A$.  

An algebra extension $A \| B$ is \textit{left depth two (D2)} if
its tensor-square $A \o_B A$ as a natural $B$-$A$-bimodule
is isomorphic to a direct summand of a finite direct sum
of the natural $B$-$A$-bimodule $A$: for some positive integer
$N$, we have 
\begin{equation}
\label{eq: D2}
A \o_B A \oplus * \cong A^N
\end{equation}
An extension $A \| B$ is \textit{right D2} if eq.~(\ref{eq: D2}) holds instead as natural
$A$-$B$-bimodules. An algebra extension is of course D2 if it is both left D2 and right D2. For example, 
if $A \| B$ is a faithfully flat algebra
(so $B$ is commutative, maps into the center of $A$
and the module $A_B$ is faithfully flat), then it 
is D2 if and only if $A$ is f.g.\ projective over
 $B$, since $A \o_B A_A$ is f.g.\ projective.  

 Since condition~(\ref{eq: D2})
implies maps in  two hom-groups satisfying $\sum_{i=1}^N g_i \circ f_i = \id_{A \o_B A}$, where  $g_i \in \Hom ({}_BA_A, {}_BA\o_B A_A) \cong (A \o_B A)^B$ (via $g \mapsto g(1)$) and $$f_i \in \Hom ({}_BA\o_B A_A, {}_BA_A) \cong
\End {}_BA_B := S $$
via $f \mapsto (a \mapsto f(a \o_B 1))$, we obtain an equivalent condition for extension $A \| B$ to be
left D2: there is a positive integer $N$, $\beta_1,\ldots,\beta_N \in S$
and $t_1,\ldots,t_N \in (A \o_B A)^B$ (i.e., satisfying for each $i = 1,\ldots,N$, $bt_i = t_i b$ for every $b \in B$) such that 
\begin{equation}
\label{eq: lqbd2}
\sum_{i=1}^N t_i \beta_i(x)y = x \o_B y
\end{equation}
for all $x,y \in A$.

To see that the quasibases equation above
is equivalent with the definition of left D2 in
eq.~(\ref{eq: D2}) it remains to note the split epi of
$B$-$A$-bimodules given by
\begin{equation}
A^N \longrightarrow A \o_B A, \ \ \ (a_1,\ldots,a_N) \longmapsto \sum_{i=1}^N t_ia_i,
\end{equation}
with splitting map $x \o y \mapsto (\beta_1(x)y,\ldots,\beta_N(x)y)$.  

Eq.~(\ref{eq: lqbd2}) is quite useful; for example, to show
$S$ finite projective as a left $V_A(B)$-module (with module action
given by $r \cdot \alpha = \lambda_r \circ \alpha$), apply 
$\alpha \in S$ to the first tensorands of the equation, set $y = 1$ and apply the
multiplication mapping $\mu: A \o_B A \to A$ to obtain
\begin{equation}
\label{eq: tys}
 \alpha(x)  = \sum_i \alpha(t^1_i) t_i^2 \beta_i(x), 
\end{equation}
where we suppress a possible summation in $t_i \in A \o_B A$ using a Sweedler notation, 
$t_i = t_i^1 \o_B t^2_i$.
But for each $i = 1, \ldots,N$, we note that $$T_i(\alpha) := \alpha(t^1_i) t^2_i \in V_A(B) := R $$ 
defines a homorphism $T_i \in \Hom ({}_RS, {}_RR)$, so that
eq.~(\ref{eq: tys}) shows that $\{ T_i \}$, $\{ \beta_i \}$ are finite dual bases for ${}_RS$.  

Similarly, an algebra extension $A \| B$ is right D2 if there is a positive
integer $N$, elements $\gamma_j \in \End {}_BA_B$ and
$u_j \in (A \o_B A)^B$ such that
\begin{equation}
\label{eq: right D2}
x \o_B y = \sum_{j=1}^N x \gamma_j(y)u_j
\end{equation}
for all $x,y \in A$.  We call the elements $\gamma_j \in S$ and $u_j \in (A \o_B A)^B$
right D2 quasibases for the extension $A \| B$. Fix this notation and the corresponding notation
$\beta_i \in S$ and $t_i \in (A \o_B A)^B$ for left D2 quasibases throughout this paper.

In the paper \cite{KK} it is shown that a subgroup $H$ of a finite group $G$ has
complex group algebras $\C H \subseteq \C G$ of depth two if and only if $H$ is normal in $G$.
From this fact we draw several examples to show that given an intermediate algebra $B \subset C \subset A$
there is in general no subalgebra pair $B \subset C$, $C \subset A$, or $B \subset A$,
where being depth two will imply another subalgebra pair in these is depth two.  For example,
the trivial subgroup is normal in a finite group $G$ containing a non-normal subgroup $H$,
so that $A \supset B$ being D2 does not imply that $A \supset C$ is D2 (unlike separability). 

If instead we ask if there may not be an exception to this rule if $C$ bears some special
relationship to $B$ within $A$, a natural candidate comes to mind as the double centralizer
of $B$ in $A$.  For  we consider the ``toy model'' for D2 extensions, a type of ``depth one'' extension called H-separable extension.  A result of Sugano and Hirata states that if $A \| B$ is H-separable,
then $A \| V := V_A(V_A(B))$ is H-separable.  The
following example shows that this does not carry
over verbatim to depth two extensions.  
   
\begin{example}
\begin{rm}
Let $A = E(W_2)$, the exterior algebra of a vector space $W_2$ (over a field $K$ of characteristic
unequal to two) with basis
$\{ e_1, e_2 \}$.  Let $B$ be the unit subalgebra, then its double centralizer
$V$ is  the center of $A$, which is $$ V = K \cdot 1_A + K \cdot e_1 \wedge e_2. $$
The extension $A \| B$ is D2 as is any finitely generated projective algebra.
However,  since $A \| V$ is a split extension, if it were D2  the module
$A_V$ would be projective, hence free since $V$ is a local algebra, which
  leads to a contradiction.  
\end{rm}
\end{example} 

We need to add a hypothesis in order to obtain some
result, such as $V \| B$ is a separable extension.  
In preparation for the next theorem, we make a definition of a weakened notion of depth two,
which is analogous to  the weakening of H-separability
in the notion of strong separability introduced in
McMahon and Mewborn \cite{MM}.  

\begin{definition}
\begin{rm}
An algebra extension $A \| B$ (with centralizer
denoted by $R$) is a \textit{weak left
depth two} extension if the module ${}_RS$ is finitely generated (f.g.) projective and the left canonical $B$-$A$-homomorphism
$\Psi: A \o_B A \to \Hom ({}_RS, {}_RA)$ defined below in
eq.~(\ref{eq: psi}) is a split epi. The definition
of a weak right D2 extension is defined oppositely.
A weak depth two extension is weak right and left D2.    
\end{rm} 
\end{definition}
A left D2 extension $A \| B$ is weak left D2 since
${}_RS$ f.g.\ projective and $\Psi$ an isomorphism
characterize left D2 extensions by \cite[2.1 (3)]{KK}. 
As an aside, the definition may be used to derive a split monomorphism $S \o_R S \into$ \newline $ \Hom ({}_BA\! \o_B\!
A_B, {}_BA_B)$ using \cite[Prop.\ 20.11]{AF} (but with
a different right $R$-module $S$ than the one used
below). The right $R$-module structure on $S$ we will use for the rest of this paper is given by $\alpha \cdot r = \alpha(-)r$ for $r \in R = V_A(B)$, $\alpha \in \End {}_BA_B = S$. 
\begin{theorem}
\label{th-closure}
Let $W : = V_A(V_A(B))$.  If $A \| B$ is a  left
(or right) D2 extension and $W \| B$
is a separable extension, then $A \| W$ is a weak left
(or right) D2 extension.  
\end{theorem}
\begin{proof}
We note that $R = V_A(B)$ and $W = V_A(R)$ are in general each other's
centralizers since also $V_A(W) = R$. Now consider the bimodule
endomorphism algebra $S'$ for the extension $A \| W$ with
base algebra $V_A(W) = R$.  We claim that the natural inclusion $\iota:\, \End {}_WA_W \into \End {}_BA_B$ is a 
split $R$-$R$-monomorphism.  Let $e = \overline{e} = e^1 \o e^2 \in W \o_B W$ be a separability element for $W \| B$.  Define a mapping $\eta:\, S \to S'$ by
$$ \eta(\alpha) := \overline{\alpha} = e^1 \alpha(e^2 ? \overline{e}^1) \overline{e}^2 $$
We note that $\eta(\lambda_r \rho_s \alpha) = \lambda_r \rho_s \eta(\alpha)$ for every $\alpha \in S$, $r,s \in R$
since elements in $R$ and $W$ commute. The mapping $\eta$ is a splitting of $\iota$ since 
for every $\beta \in S'$ we have $\overline{\beta} = \beta$
as $\beta$ is $W$-linear and $e^1 e^2 = 1$. 

We will now show that $A \| B$ left D2 implies $A \| W$ is weak left
D2.  Since $A \| B$ is left D2, the module ${}_RS$ is finite projective.  Since ${}_RS'$ is a direct summand of ${}_RS$, it
too is finite projective.  It will then suffice to show
that the mapping 
\begin{equation}
\label{eq: psi}
\Psi: A \o_W A \to \Hom ({}_RS', {}_RA),
\ \ \ \Psi(x \o y)(\beta) := \beta(x)y
\end{equation}
is a split $W$-$A$-epimorphism. 

Define a splitting map $$\Hom ({}_RS', {}_RA) \to A \o_W A, \ \ \ 
 G \longmapsto \sum_i e^1 t^1_i \o_W t^2_i e^2 G(\overline{\beta_i}) .$$ 
Indeed, for each $\beta \in \End {}_WA_W$, 
\begin{eqnarray*}
\Psi(\sum_i e^1 t^1_i \o_W t^2_i e^2 G(\overline{\beta_i}))(\beta)
& = & \sum_i \beta(e^1 t^1_i) t^2_i e^2 G(\overline{\beta_i}) \\
& = & \sum_i G(e^1 \beta(t^1_i)t^2_i\beta_i(e^2 ? \overline{e}^1) \overline{e}^2) \\
& = & G(e^1 \beta(e^2 ? \overline{e}^1) \overline{e}^2) = G(\beta)
\end{eqnarray*}
since $\beta(t^1_i)t^2_i \in R$ for each $i$, $G$ is left $R$-linear and $\sum_i \beta(t^1_i) t^2_i \beta_i(x) = \beta(x)$ for each $x \in A$ by eq.~(\ref{eq: tys}).  

The proof that  $W \| B$ separable and $A \| B$  right D2 implies  that $A \| W$ is weak right D2, is  similar.  
\end{proof}

For example, the theorem is well-known for the complex group algebras
corresponding to the subgroup situation
where $H$  normal in a finite group $G$ implies that
its centralizer $V_G(H)$ is normal in $G$.  
We provide a characterization for a weak left D2 extension.  Note that the $B$-$A$-subbimodule $U$ below coincides with the reject of $A$ in $A \o_B A$
\cite[p.\ 109]{AF}.

\begin{prop}
An algebra extension $A \| B$ is weak left D2 if and only if 
as $B$-$A$-bimodules, the three conditions below are satisfied: 
\begin{enumerate}
\item $ A \o_B A = U \oplus L  $, where 
\item  $\Hom (U, A) = \{ 0 \} $,  and 
\item $L \oplus * \cong A^N$ for some positive integer $N$.
\end{enumerate}
\end{prop}
\begin{proof}
($\Rightarrow$)  Since ${}_RS$ is f.g.\ projective, there is a positive integer $N$
such that ${}_RS \oplus * \cong {}_RR^N$.  Then $\Hom ({}_RS, {}_RA) \oplus * \cong A^N$
as $B$-$A$-bimodules.  But $\Psi: A \o_B A \to \Hom ({}_RS, {}_RA)$ defined
by $\Psi(x \o y)(\alpha) = \alpha(x)y$ is a split epi.  Let $L$ be the image  in $A \o_B A$
of $\Hom ({}_RS, {}_RA)$ under a split monic, therefore
satisfying condition (3). Let $U = \ker \Psi$.  Then
$A\o_B A = U \oplus L$.  We next show that ${\rm Hom}_{B-A}\, (U,A) = 0$
where the unlabelled homomorphism groups are w.r.t.\ $B$-$A$-bimodules:
$$ S \stackrel{\cong}{\longrightarrow} \Hom (A \o_B A, A) \cong  
\Hom (\Hom ({}_RS, {}_RA),A) \oplus \Hom (U , A) \cong S \oplus \Hom (U,A) $$
since ${}_RS$ f.g.\ projective implies $\Hom (\Hom ({}_RS, {}_RA),A) \cong \Hom (A,A) \o_RS \cong
S$ by a standard isomorphism (e.g., \cite[3.4]{MM}) and the fact that $\End ({}_BA_A) \cong R$
via $f \mapsto f(1)$. The first arrow above is given by $\alpha \mapsto (x \o_B y \mapsto \alpha(x)y)$ with inverse $F \mapsto F(? \o 1_A)$.  It is not hard to check that the composite
isomorphism is the identity on the direct summand $S$, whence $\Hom (U, A) = 0$.  

($\Leftarrow$) Since $S \cong \Hom (A \o_B A, A)$, it follows from substitution of condition
(1)  and applying condition (2) that $S \cong \Hom ({}_BL_A,{}_BA_A)$.  From condition (3)
and a derivation as in that of eq.~(\ref{eq: lqbd2}), we arrive at elements $2N$ elements 
$v_j = v_j^1 \o v_j^2 \in L^B \subseteq (A \o_B A)^B$, $\delta_j \in S$ such that
\begin{equation}
\label{eq: partial qb}
\ell = \sum_{j=1}^N v^1_j \o_B v^2_j \delta_j(\ell^1)\ell^2
\end{equation}
for all $\ell = \ell^1 \o \ell^2 \in L \subseteq A \o_B A$.  Again from condition (3), $\Hom (L,A) \oplus * \cong \End ({}_BA_A)^N$,
whence ${}_RS \oplus * \cong {}_RR^N$ and $S$ is left f.g.\ projective $R$-module.  

The mapping $\Psi: A \o_B A \to \Hom ({}_RS, {}_RA)$ is split by the $B$-$A$-bimodule
homomorphism 
$$\sigma: \Hom ({}_RS, {}_RA) \to A \o_B A, \ \ \ \sigma(G) := \sum_{i=1}^N v_j^1 \o_B v_j^2 G(\delta_j) $$
since $\sum_j \alpha(v^1_j)v^2_j G(\delta_j) = G(\sum_j \alpha(v^1_j)v^2_j\delta_j) = G(\alpha)$
(for all $\alpha \in S \cong \Hom (L,A)$), which follows from eq.~(\ref{eq: partial qb}).  
\end{proof}

From the proof it is clear that another characterization of weak left D2 extension $A \| B$
is that its tensor-square has a direct sum decomposition as in conditions (1) and (2),
where all elements of $L$ satisfy a left quasi-bases equation~(\ref{eq: partial qb}). The extent to which
a weak depth two extension has a Galois theory might be
an interesting problem.  

%%%%%%%%%%%%%%%%%%%%%%%%%%%%%%%%%%%%%%%%%%%%%%%%%%%%%%%%%%%%%%%%%%%%%%%%%%%%%%%%%%%%%%%%%%%%%%%%%
\section{Further closure properties of depth two with counterexamples}

Unlike separable extensions and Frobenius extensions, depth two is not a transitive property.
If $G$ is a finite group with normal subgroup $N$ having a normal
subgroup $K$ where $K \not \ract G$, then the corresponding
complex group algebras $A = \C G$, $B = \C K$, and $C = \C N$
satisfy $A \supset C \supset B$ with $A \| C$ and $C \| B$ both
D2 but $A \| B$ not D2. However,  normality of subgroups
$G \geq N \geq K$ satisfies $K \ract G$  normal $\Rightarrow$
$K \ract N$.  Thus 
it may come as a surprise that $A \| C$ and $A \| B$ D2   $\not \Rightarrow$ $C \| B$ D2, which may be seen from the example
$A= M_2(\C)$, $B = \C \times \C$ and $C = T_2(\C)$,
the triangular and full $2 \times 2$ matrix algebras and the algebra of diagonal
matrices, as shown in the next proposition.

\begin{prop}
The matrix algebras $A = M_n(k)$, $B = {\rm Diag}_n(k)$ and $C = T_n(k)$ 
over any field $k$ satisfy:  $A \| B$  and $A \| C$ are H-separable
(and therefore D2), but $C \| B$ is not D2.  
\end{prop}
\begin{proof}
An algebra extension $A \| C$ is H-separable iff
$1 \o_B 1$ may be expressed as a sum of products of elements in the centralizer $V_A(C)$ and (Casimir)
elements in $(A \o_C A)^A$. Recall that for any fixed $j = 1, \ldots,$ or $n$,
$\sum_{i = 1}^n e_{ij} \o_k e_{ji}$ is a Casimir element in $(A \o_k A)^A$. 
But note 
that $1 \o_C 1$ reduces to a canonical image of this Casimir element, since the matrix units $e_{ii}, e_{1i} \in C$ for $1 \leq i  \leq n$  yields
$$ 1 \o_C 1 = \sum_{i=1}^n e_{ii} \o_C e_{ii} = \sum_{i=1}^n
e_{i1} \o_C e_{1i}.$$ 

For a similar reason $A \| B$ is H-separable, since the centralizer $V_A(B) = B$,
and $$ 1 \o_ B 1 = \sum_{i=1}^n e_{ii} \o_B e_{ii} =  \sum_{i=1}^n e_{ii} (\sum_{j=1}^n e_{ji} \o_B e_{ij}).$$

Now, if an extension $C \| B$ is D2, then its centralizer $R := V_C(B)$ is a normal subalgebra in $C$:
i.e., for each two-sided ideal $I$ in $C$, we have
the $C$-invariance of the  ideal contracted to $R$,
$(R \cap I)C = C( R \cap I)$ \cite[Prop.\ 4.2]{LK2005}.
But $R = B$ in our example and $C$  has the two-sided ideal $I = \sum_{i=1}^n ke_{1i}$ where $C(I \cap B) \neq (I \cap B)C$. Hence, $C \| B$ is not D2.  
\end{proof}

It is similarly shown that $C \| B$ is not one-sided D2.  
There is a certain transitivity of the depth two property
when it follows an H-separable extension.

\begin{prop}
If the algebra extension $A \| C$ is right (or left) D2 and 
the extension $C \| B$ is H-separable, then $A \| B$ is right
(or left) D2.
\end{prop}
\begin{proof}
If $C \| B$ is H-separable, we have $${}_CC \o_B C_C \oplus * \cong
{}_C{C^N}_C$$ for some positive integer $N$.  Apply the
functor ${}_C\mathcal{M}_C \to {}_A\mathcal{M}_C$ given
by \newline
 ${}_AA \o_C \, ? \, \o_C A_C$, to obtain
\begin{equation}
\label{eq: tys2}
{}_A A \o_B A_C \oplus * \cong \oplus^N {}_AA \o_C A_C. 
\end{equation}
If $A \| C$ is right D2, we have ${}_A A \o_C A_C \oplus * \cong
{}_A {A^M}_C$ for some positive integer $M$.  It follows from this applied to eq.~(\ref{eq: tys2}), then restricting from right $C$-modules
to $B$-modules that
$$ {}_A A \o_B A_B \oplus * \cong {}_A{A^{NM}}_B, $$
which is the condition that $A \| B$ is right D2.  It is similarly
proven that a left
D2 following  an H-separable extension is altogether left D2.
\end{proof}
There has been a question of whether a
 right or left progenerator H-separable extension $A \| B$ is split 
(i.e., has a $B$-$B$-bimodule projection $A \to B$), whence Frobenius: an affirmative answer implies some generalizations of results of Noether-Brauer-Artin on simple algebras \cite{Su99}. Unfortunately, the next example, derived from  the endomorphism ring theorem for D2 extensions in \cite{LK2006},
of a one-sided free H-separable non-Frobenius extension
rules out this possibility.      
  
\begin{example}
\begin{rm}
Let $K$ be a field and $B$ the $3$-dimensional algebra of  upper triangular $2 \times 2$-matrices, which is not self-injective.
Since $B \| K1$ is trivially D2, the endomorphism algebra $A := \End B_K \cong M_3(K)$ is a left
D2 extension of $\lambda(B)$, which w.r.t.\ the ordered basis $\bra e_{11}, e_{12}, e_{22} \ket$
$\lambda(B)$ is the subalgebra of matrices
$$ [x,y,z] :=  \left( \begin{array}{ccc}
x & 0 & 0 \\
0 & x & y \\
0 & 0 & z
\end{array}
\right)
$$  
Now the centralizer $R$ is the $3$-dimensional algebra spanned by matrix units $e_{11}$, $e_{21}$, $e_{22} + e_{33}$.  
The module ${}_BA$ is free with left $B$-module isomorphism $A \to B^3$ given by ``separating out the columns''
$$(a_{ij}) \mapsto
([a_{11},a_{21},a_{31}], [a_{12},a_{22}, a_{32}], [a_{13}, a_{23}, a_{33}])$$
whence $A \o_B A$ and $\Hom (R_K,A_K)$ are both $27$-dimensional.  The $A$-$A$-homomorphism 
$A \o_B A \to \Hom (R_K, A_K)$ given by $a \o c \mapsto (r \mapsto arc)$ is easily computed to 
be surjective, therefore an isomorphism, whence ${}_A A \o_B A_A \cong {}_AA^3_A$, which shows
$A \| B$ is H-separable (and D2). The extension $A \| B$ is not Frobenius
since $B$ is not a Frobenius algebra; therefore $A \| B$ is not split. (Alternatively, if
 there is a $B$-linear projection $E: A \to B$, we note $E(e_{32}) = 0$, so $e_{33} = e_{32} e_{23} \in \ker E$,
a contradiction.) 
 By applying the matrix transpose, the results of this example may be transposed to a right-sided version. 
\end{rm}
\end{example} 
  
Next we deal with elementary ways to generate new depth
two extensions from old ones.  In our first case, we note an extension is depth two if and only if all its components in a finite direct product are D2.  
 
\begin{prop}
Suppose $B_k \rightarrow A_k$ is an algebra homomorphism
$\iota_k$ for each $k = 1, \ldots, n$. Let $B = B_1 \times \cdots \times B_n$, $A = A_1 \times \cdots A_n$ and $B \rightarrow A$ be induced by $\iota= \times_{i=1}^n \iota_i$.  Then $A \| B$ is D2 if and only if $A_k \| B_k$ is D2 for each $k = 1, \ldots,n$.
\end{prop}

\begin{proof}
Let $p_i: A \to A_i$ and $\sigma_i: A_i \to A$ be
the canonical algebra morphisms satisfying $ \sum_{i=1}^n \sigma_i \circ p_i = \id_A$, $p_i \circ \sigma_j =  \id_{A_i}\delta_{ij}$.  Let $e_i = \sigma_i \circ p_i(1_A)$ be the canonical orthogonal
idempotents.  Similarly, let $\pi_i: B \to B_i$ and
$\eta_i: B_i \to B$ be the corresponding canonical mappings
and orthogonal idempotents $f_i = \eta_i \circ \pi_i(1_B)$.  These satisfy commutative squares corresponding to $p_i \circ \iota = \iota_i \circ  \pi_i$, $\sigma_i \circ \iota_i = \iota \circ \eta_i$ and $\iota(f_i) = e_i$ 
for each $i = 1,\ldots,n$.  

Note that $p_i \o p_i$ induces a $B$-$B$-bimodule
epimorphism from $A \o_B A \to A_i \o_{B_i} A_i$
split by $\sigma_i \o \sigma_i$.  Since $\sigma_i(A_i) \o_B \sigma_j(A_j) = 0$ where $i \neq j$ by using $e_i$, we see that $(A \o_B A)^B \cong \oplus_{i=1}^n (A_i \o_{B_i}  A_i)^{B_i}$. 

Similarly, $\End {}_BA_B \cong \times_{i=1}^n \End {}_{B_i}(A_i)_{B_i}$ via $\alpha \mapsto (p_i \circ \alpha \circ \sigma_i)_{i=1}^n$. Also, the centralizer $V_A(B) = \times_{i=1}^n V_{A_i}(B_i)$. 

 The proof is now completed
by projecting eq.~(\ref{eq: lqbd2}) via $p_i \o p_i$
onto $n$ left D2 quasibases equations for $A_i \| B_i$, 
and conversely gluing together $n$ left D2 quasibases equations into one for $A \| B$.  A similar argument
using the right D2 quasibases eq.~(\ref{eq: right D2})
shows $A \| B$ is right D2 $\Leftrightarrow$ each
$A_i \| B_i$ is right D2.  
\end{proof}

The proof shows that the $R$-bialgebroids $S$ and $T$
for $A \| B$ are not surprisingly  direct products of $R_i$-bialgebroids
$S_i = \End {}_{B_i}(A_i)_{B_i}$ and $T_i =(A_i \o_{B_i} A_i)^{B_i}$, respectively, where $R_i = V_{A_i}(B_i)$.
We retain this notation in dealing  with the tensor product of algebras next.  Tensor
product of finitely many D2 algebra extensions is D2, but the converse is more demanding and requires several hypotheses.
Suppose then that all the algebras in the proposition
below are finite dimensional algebras over a field $K$ and a reminder that unadorned tensors are over $K$. 

\begin{prop}
Let $B_k \rightarrow A_k$ be algebra homomorphisms
$\iota_k$ for each $k = 1, \ldots, n$ where
each $B_k$ is a separable algebra. Let $B = B_1 \otimes \cdots \o B_n$, $A = A_1 \o \cdots \o A_n$ and $B \rightarrow A$ be induced by $\iota= \o_{i=1}^n \iota_i$.  Then $A \| B$ is D2 if and only if $A_k \| B_k$ is D2 for each $k = 1, \ldots,n$.
  \end{prop}

\begin{proof}
($\Leftarrow$)  Here no hypotheses are required beyond
the objects defined.  We note the simple rearrangement mapping  $\otimes_{i=1}^n T_i \to T$, as well as $\otimes_{i=1}^n S_i \to S$ given by sending $\alpha_1 \o \cdots
\o \alpha_n$ to itself (for $\alpha_i \in S_i$). Let $u_{ij} \in T_i$ and $\gamma_i \in S_i$ be given
for each $i$ satisfying the right D2 quasibases eq.\
$$ a_i \o_{B_i} {a_i}' = \sum_{j=1}^n a_i \gamma_{ij}({a_i}') u_{ij}^1 \o_{B_i} u_{ij}^2, $$
for $a_i, {a_i}' \in A_i$.  Then $\gamma_j = \otimes_{i=1}^n \gamma_{ij} \in S$ and $$u_j = (u_{1j}^1 \o_K \cdots \o_K u_{nj}^1) \o_B (u^2_{1j} \o_K \cdots \o_K u^2_{nj}) \in T $$
as well as $a = a_1 \o_K \cdots \o_K a_n$,
$a' = {a_1}' \o_K \cdots \o_K {a_n}'$ satisfy
$$ a \o_B a' = \sum_j a \gamma_j(a') u^1_j \o_B u_j^2.$$
A similar argument shows $A \| B$ is also left D2.  

($\Rightarrow$) We need a lemma \cite[2.4]{DI} stating
that for $A$ and $B$ two $K$-algebras with finite projective modules ${}_AM$ and ${}_BN$ and two
others ${}_AM'$ and ${}_BN'$, the natural
mapping 
\begin{equation}
\label{eq: DeMeyer-Ingraham}
\Hom ({}_AM, {}_AM') \o_K \Hom ({}_BN, {}_BN')
\stackrel{\cong}{\longrightarrow} \mbox{Hom}_{A \o B} (M\o N, M' \o N')
\end{equation}
is an isomorphism (with inverse $F \mapsto \sum_{i,j}
f_i \o g_j F(m_i \o n_i)$ where $m_i, f_i$ are dual
bases for ${}_AM$ and $n_j, g_j$ are dual bases for
${}_BN$). This fact may be extended by induction to any
number of tensor factors.  

Since each $B_i^e := B_i \o B_i^{\rm op}$ is semisimple,
the $B_i^e$-modules $A_i$ are finite projective, whence
$$S_1 \o_K \cdots \o_K S_n \cong S, $$
since $B^e \cong B_1^e \o \cdots \o B_n^e$.  

Next recall that for any algebra $B$ and every $B$-$B$-bimodule $V$
the  subgroup of $B$-central elements, $V^B \cong \mbox{\rm Hom}_{B^e}(B,V)$ in a functorial way,
which is right exact if $B$ is $B^e$-projective, i.e.,
$B$ is a separable algebra \cite{DI}.  It then follows
from $A \o_B A \cong \o_{i=1}^n A_i \o_{B_i} A_i$ that 
$$T \cong \mbox{\rm Hom}_{B^e}(B, A \o_B A) \cong
\otimes_{i=1}^n \mbox{\rm Hom}_{B_i^e}(B_i, A_i \o_{B_i} A_i) = T_1 \o \cdots \o T_n. $$
Similarly, $R = A^B = R_1 \o \cdots \o R_n$.

We will use the characterization that $A \| B$ is left
D2 if $T_R$ is finite projective and $T \o_R A \stackrel{\cong}{\longrightarrow} A \o_B A$ via the map
$\beta(t \o_R a) := t^1 \o_B t^2a$ \cite[2.1(4)]{KK}. 
Under the decompositions above, $\beta$ decomposes
into corresponding mappings $$\o_{i=1}^n \beta_i:
\o_{i=1}^n T_i \o_{R_i} A_i \stackrel{\cong}{\longrightarrow} \o_{i=1}^n A_i \o_{B_i} A_i. $$  
By finite dimensionality one shows that each
$\beta_i$ is injective and surjective. Finally,
$T_R$ is finite projective, so $T_{R_i}$ is finite
projective, whence the summand $T_i$ is finite projective as a right $R_i$-module.  Thus
$A_i \| B_i$ is left D2 for each $i = 1, \ldots,n$,
and by a similar argument, it is right D2.  
  \end{proof}

From the proof we note that the $R$-bialgebroids $S$
and $T$ for a (one-sided or two-sided) D2 extension $A \| B$ once again decompose, this time into a  tensor
product of $R_i$-bialgebroids $S_i$ and $T_i$
(with antipode if the characteristic of $K$ is zero \cite[3.6]{LK2006}.
Note that any algebra is a D2 extension of itself.  
\begin{cor}
Let $B$ be a separable algebra. Then 
$A \| B$ is D2 $\Leftrightarrow$ $M_n(A) \| M_n(B)$
is D2. 
\end{cor}
 
%%%%%%%%%%%%%%%%%%%%%%%%%%%%%%%%%%%%%%%%%%%%%%%%%%%%%%%%%%%%%%%%%%%%%%%%%%%%%%%%%%%%%%%%%%%%%%%%%%%%%%%%%%%%%%%

\section{A characterization of Galois extension}

Tachikawa  \cite{T} studies the double centralizer condition for modules and the property
of being balanced for rings.
We give a related result --- that a subalgebra satisfying the  double centralizer (or bicommutant)
condition has a balanced right or left regular representation on the over-algebra. 

\begin{lemma}
\label{lem-bekvem}
If $B \subseteq A$ is a subalgebra, then  
$$B \ \subseteq \ A^S := \{ x \in A \, : \, \forall \alpha \in \End {}_BA_B, \, \alpha(x) = \alpha(1_A)x \} \subseteq \ V_A(V_A(B)). $$
If $B$ satisfies the double
centralizer  condition, $V_A(V_A(B)) = A$,
or equals the invariant subalgebra $A^S = B$, 
then the natural modules $A_B$ and ${}_BA$ are balanced. 
\end{lemma} 
\begin{proof}
Again let $S$ denote $\End {}_BA_B$.  Clearly $B \subseteq A^S$.  Assume $x \in A^S$.
Then for each $r \in V_A(B)$ we have $\rho_r \in S$, so $\rho_r(x) = \rho_r(1) x$,
i.e. $xr = rx$, so $x \in V_A(V_A(B))$.  Hence,
$B \subseteq A^S \subseteq V_A(V_A(B))$.  

Suppose $B$ satisfies the double centralizer condition
in $A$.  Then $A^S = B$.  
To see that $A_B$ is balanced, let $E$ denote $\End A_B$ and consider $f \in \End {}_EA$. Then $\forall a \in A$,
$$ f(a) = f(\lambda_a(1)) = \lambda_a(f(1)) = a f(1), $$
thus $x := f(1)$ satisfies $f = \rho_x$.  It suffices to  show that $x \in B$. Let $\alpha \in S$, then
$$ \alpha(x) = f(\alpha(1)) = \lambda_{\alpha(1)} f(1) = \alpha(1)x$$
whence $x \in A^S = B$. The proof that ${}_BA$ is balanced is similar.   
\end{proof} 

Recall that an algebra extension $A \| B$ is a \textit{$\mathcal{G}$-Galois extension} \cite{Su} if $\mathcal{G}$ is
a group of automorphisms of the algebra $A$ fixing
each element in $B$ such that
\begin{enumerate}
\item the natural module $A_B$ is finitely generated
and projective;
\item the mapping $j: A \rtimes G \to \End A_B$
given by $j(a \# \sigma)(x) = a \sigma(x)$ (for each 
$a,x \in A$, $\sigma \in G$) is an isomorphism;
\item the set of invariants
$A^G = \{ x \in A \, : \, \sigma(x) = x, \ \forall\,
\sigma \in G \} $ is equal to $B$: $A^G = B$.  
 \end{enumerate} 

Now consider left or right Galois extensions for bialgebroids
and their characterization as left or right depth two
and balanced extensions \cite{fer}.  
Next we give a characterization for finite projective
Galois extensions which is very similar to the one above
for group-Galois extensions:

\begin{theorem}
Suppose $A \| B$ is an algebra extension with centralizer denoted by 
$R = V_A(B)$ and $A_B$ finite projective. Then $A \| B$ is a left 
Galois extension iff ${}_RS$ is finite projective,
$j: A \o_R S \to \End A_B$ given
by $j(a \o \alpha)(x) = a\alpha(x)$ is an isomorphism,
and $A^S = B$.
\end{theorem}
\begin{proof}
($\Rightarrow$) From \cite[Theorem 2.1]{fer},
$A \| B$ is left D2 and balanced, and from \cite[3.10, 4.1]{KS}, $j$ is an isomorphism and
$A^S = B$.  That ${}_RS$ is finite projective is noted
above after eq.~(\ref{eq: tys}).

($\Leftarrow$) If ${}_RS \oplus * \cong {}_RR^N$,
we tensor by $A \o_R -$ and apply the $A$-$B$-isomorphism $j$ to obtain $\End A_B \oplus * \cong A^N$ as natural $A$-$B$-bimodules.  Since $A_B$
is finite projective, we may apply \cite[3.8]{LK2005} to see that $A \| B$ is left D2.  Since $A^S = B$,  Lemma~\ref{lem-bekvem} informs us that $A$ is balanced
over $B$.  Thus, $A \| B$ is a left Galois extension by \cite[Theorem 2.1]{fer}.
\end{proof}

%%%%%%%%%%%%%%%%%%%%%%%%%%%%%%%%%%%%%%%%%%%%%%%%%%%%%%%%%%%%%%%%%%%%%%%%%%%%%%%%%%%%%%%%%%%%%%%

\section{Depth two and cup product in simplicial Hochschild cohomology}

Let $A \| B$ be an extension of  $K$-algebras.  We briefly recall the $B$-relative
Hochschild cohomology of $A$ with coefficients in $A$ (for coefficients in a bimodule,
see the source \cite{H}).  The zero'th cochain group $C^0(A,B;A) =
A^B = R$, while the $n$'th cochain group $C^n(A,B;A) = \mbox{\rm Hom}_{B-B}(A \o_B \cdots \o_B A,A)$ ($n$ times $A$ in the domain).  In particular, $C^1(A,B;A) = S= \End {}_BA_B$.  The
coboundary $\delta_n: C^n(A,B;A) \to C^{n+1}(A,B;A)$ is given by
\begin{eqnarray}
\label{eq: Hoch}
(\delta_n f)(a_1 \o \cdots \o a_{n+1}) &=& a_1 f(a_2 \o \cdots \o a_{n+1}) 
+ (-1)^{n+1} f(a_1 \o \cdots \o a_n)a_{n+1}   \nonumber \\
 & + & \sum_{i=1}^n  
(-1)^i f(a_1 \o \cdots \o a_i a_{i+1} \o \cdots \o a_{n+1})  
\end{eqnarray}
which satisfies $\delta_{n+1} \circ \delta_n = 0$ for each positive integer $n $.  Its cohomology is denoted by
$H^n(A,B;A) = \ker \delta_n / \Im \delta_{n-1}$, and might be referred to as a 
simplicial Hochschild cohomology, since  this cohomology
is isomorphic to simplicial cohomology
if $A$ is the poset algebra of a simplicial complex \cite{GS}. 

The cup product $\cup: C^m(A,B;A) \o_K C^n(A,B;A) \to C^{n+m}(A,B;A)$ makes use of the multiplicative
stucture on $A$ and is given by
\begin{equation}
(f \cup g)(a_1 \o \cdots \o a_{n+m}) = f(a\1 \o \cdots \o a_m)g(a_{m+1} \o \cdots \o a_{n+m})
\end{equation}
which satisfies the equation $\delta_{n+m}(f \cup g) = (\delta_m f) \cup g + (-1)^m f \cup \delta_n g$ \cite{GS}.  Cup product therefore passes to a product on the cohomology.  
  We note that $(C^*(A,B;A), \cup, +, \delta)$ is a differential graded algebra (perhaps
negatively graded according to the convention used)
which we denote by $D(A,B)$.  

\begin{theorem}
Suppose $A \| B$ is a right or left D2 algebra extension.  Then $D(A,B)$ is generated as an
algebra by its degree one elements and is isomorphic to the tensor algebra on $C^1(A,B;A)$
over $C^0(A,B;A)$.  
\end{theorem}

\begin{proof}
The idea of the proof is to generalize the isomorphism $$S \o_R S \stackrel{\cong}{\longrightarrow} \Hom ({}_BA \o_B A_B, {}_BA_B)$$
via $\alpha_1 \o_R \alpha_2 \mapsto \alpha_1 \cup \alpha_2$ \cite[3.11]{KS} in the notation
for $S = C^1(A,B;A)$ and $R = C^0(A,B;A)$ above.  This shows that any $2$-cochain is the cup
product of 1-cochains.  Similarly, any $n$-cochain is the cup product of $1$-cochains, since
\begin{equation}
S \o_R \cdots \o_R S \stackrel{\cong}{\longrightarrow} \mbox{\rm Hom}_{B-B}(A \o_B \cdots \o_B A, A).
\end{equation}
 via $\alpha_1 \o \cdots \o \alpha_n \mapsto \alpha_1 \cup \cdots \cup \alpha_n$.
We prove this by induction on $n$, the statement  holding in $n= 1,2$.  Suppose it holds
for $n < m$.  Note that $S \o_R C^{m-1}(A,B;A) \stackrel{\cong}{\longrightarrow} C^m(A,B;A)$
via $\alpha \o g \mapsto \alpha \cup g$ since an inverse is clearly given by
$f \longmapsto$ \newline $ \sum_j \gamma_j \o_R u^1_j f(u^2_j \o_B \cdots \o_B-)$ in terms of right D2 quasibases.
By the induction hypothesis $C^{m-1}(A,B;A) \cong S \o_R \cdots \o_R S$ ($m-1$ times $S$) via the cup product, so the proof is complete.  It follows that there is an isomorphism of algebras, 
\begin{equation}
T^*({}_RS_R) \stackrel{\cong}{\longrightarrow} D(A,B)
\end{equation}
via the cup product mapping above.  
\end{proof} 

The inverse mapping for  $S \o_R S \o_R S \cong \Hom ({}_BA \o_B A \o_B A_B, {}_BA_B)$
implied by the proof has a different expression than that in
\cite[(46)]{fer}: it comes out here as
$g \mapsto \sum_{i,j} \gamma_j \o_R \gamma_i \o_R u^1_i u_j^1 g(u^2_j \o_B u^2_i \o_B -)$
for $g \in C^3(A,B;A)$.  
%%%%%%%%%%%%%%%%%%%%%%%%%%%%%%%%%%%%%%%%%%%%%%%%%%%%%%%%%%%%%%%%%%%%%%%%%%%%%%%%%%%%%%%%%%%%%%%%%%%%
\section{Codepth two coalgebra homomorphisms and their bialgebroids}

In this section, we dualize the notion of depth two for algebra homomorphisms to obtain a notion of codepth two for coalgebra homomorphisms.  We obtain workable codepth two quasibases via simplifications of certain  hom-groups of comodule homomorphism.  We then establish
 a right bialgebroid structure on the bicomodule endomorphisms, where the base algebra is the centralizer of the dual algebra homomorphism.  

Let $g: C \to D$ be a homomorphism of coalgebras over a field (alternatively, coalgebras which are flat over a base ring $K$ with $K$-duals that separate points).  Then $C$ has an induced $D$-$D$-bicomodule structure given
by left coaction $$ \rho^L: C \to D \o C,\ \rho^L(c) = c\-1 \o c\0 := g(c\1) \o c\2,$$
and by right coaction $$ \rho^R: C \to C \o D,\ \rho^R(c) = c\0 \o c\1 := c\1 \o g(c\2). $$ 
In a similar way, any $C$-comodule becomes a $D$-comodule via the homomorphism $g$, 
the functor of corestriction \cite[11.9]{BW}.  We denote a right $D$-comodule $M$ by
$M^D$ for example.   This should not be confused with our notation $A^B$ for the 
centralizer of a subring $B$ in a ring $A$, which we need to work with simultaneously below.  
(Unadorned tensors between modules are over $K$, we use a generalized Sweedler notation,
the identity is sometimes denoted by its object,  
and basic terminology such as coalgebra homomorphism, comodule or bicomodule is defined in the standard way such as in \cite{BW}.) 

Recall that the cotensor product $$ C \Box_D C = \{  c \o {c'} \in C \o C \, | \,
c\1 \o g(c\2) \o c' =  c \o g({c'}\1) \o {c'}\2 \},$$
where we suppress a possible summation $c \o c' = \sum_i c_i \o {c'}_i$.  
For example, if $g = \eps: C \to K$ the counit on $C$, $C \b_D C= C \o C$.  

Note that $C \Box_D C$ is a natural $C$-$C$-bicomodule via the coproduct $\cop$ on $C$
applied as $\cop \o C$ for the left coaction and $C \o \cop$ for the right coaction \cite[11.3]{BW}.  
Then $\underline{\cop}:  \, C \to C\b_D C$ induced by $\cop$ (where $\underline{\cop}(c) := c\1 \o c\2$)
 is a $C$-$C$-bicomodule monomorphism.  As $D$-$C$-bicomodule it is split
by $c \o c' \mapsto \eps(c) c'$,
and as a $C$-$D$-bicomodule $\underline{\cop}$ is split by $c \o c' \mapsto c \eps(c')$
for $c \o c' \in C \b_D C$. (Since $c\1 \o g(c\2) \o c' = c \o g({c'}\1) \o {c'}\2$,
it follows that $g(c) \o c' = \eps(c) g({c'}\1) \o {c'}\2$,
whence $c \o c' \mapsto \eps(c)c'$ is left $D$-colinear.) 
For example, if $D = C$ and $g = \id_C$, then $C \b_C C \cong C$, since
$\underline{\cop}$ is surjective.     

It follows that $C$ is in general isomorphic to a direct summand of $C \b_D C$
as $D$-$C$-bicomodules: ${}^D C \b_D C^C \cong {}^DC^C \oplus * $. 
Left codepth two coalgebra homomorphisms have the special complementary property: 

\begin{definition}
\begin{rm}
A coalgebra homomorphism $g: C \to D$ is said to be \textit{left codepth two (cD2)}
if for some positive integer $N$,
\begin{equation}
{}^D C \b_D C^C \oplus * \cong {}^D(C^N)^C, 
\end{equation}
i.e., the cotensor product $C \b_D C$ is isomorphic to a direct summand of
a finite direct sum of $C$ with itself as $D$-$C$-bicomodules.
Right codepth two coalgebra homomorphisms are similarly defined.
\end{rm}
\end{definition}

The definition implies that there are $D$-$C$-bicomodule homomorphisms
$f_i \in {\rm Hom}^{D-C} (C, C \b_D C)$ and $g_i \in {\rm Hom}^{D-C} (C \b_D C, C)$
such that 
\begin{equation}
\label{eq: basic}
C \b_D C = \sum_{i=1}^N f_i \circ g_i.
\end{equation}  
We are then interested in obtaining simplications for these two hom-groups.  
\begin{prop}
Given a coalgebra homomorphism $g : C \to D$, 
$$ \End {}^DC^D \cong {\rm Hom}^{D-C}(C, C \b_D C). $$
\end{prop}
\begin{proof}
Define a mapping $\End {}^DC^D \to {\rm Hom}^{D-C}(C, C \b_D C)$ by
\begin{equation}
\label{eq: map1}
\alpha \mapsto (\alpha \o C)\circ \underline{\cop}.
\end{equation}
Note that $\alpha(c\1) \o c\2 \in C \b_D C$ for every $c \in C$,
since $\alpha(c\1) \o g(c\2) = $
$ (C \o g) \cop(\alpha(c))$ by right $D$-colinearity
of $\alpha$. 

This mapping has inverse ${\rm Hom}^{D-C}(C, C \b_D C) \to \End {}^DC^D$ given
by $F \mapsto (C \o \eps) \circ F$.  This is clearly left $D$-colinear, and
right $D$-colinear since $F$ and $C \o \eps$ are so.  
It is an inverse since 
$$ (C \o \eps) \circ (\alpha \o C) \circ \cop = (\alpha \o \eps) \circ \cop = \alpha $$
for each $\alpha \in \End {}^DC^D$, and right $C$-colinearity for $F \in {\rm Hom}^{D-C}(C, C \b_D C)$
means $(C \o \cop) \circ F = (F \o C)\circ \cop$, so  
$$ (((C \o \eps) \circ F) \o C) \circ \cop = (C \o \eps \o C)\circ (C \o \cop) \circ F = F. \qed $$
\renewcommand{\qed}{}\end{proof}

Part of this proposition may be derived directly from
the hom-cotensor relation \cite[11.10]{BW}.  

Let $D^*$ and $C^*$ be the dual $K$-algebras of coalgebras $D$ and $C$, respectively,
with multiplication given by the convolution product and unity element equal to the counit.
Note that $g: C \to D$ induces the algebra homomorphism $g^*: D^* \to C^*$
given by $d^* \mapsto d^* \circ g$.  Then $C^*$ obtains a $D^*$-$D^*$-bimodule structure via
$g^*$ and we define the centralizer $V_{C^*}(D^*)$ to be the set of all elements $c^* \in C^*$ such
that for all $d^* \in D^*$, $c^* g^*( d^*) = g^*(d^*) c^*$ or equivalently 
suppressing $g^*$, $c^* \cdot d^* = d^* \cdot c^*$. 

\begin{lemma}
Given coalgebra homomorphism $g: C \to D$.  Then $\End {}^DC^C \cong V_{C^*}(D^*)$
via $f \mapsto  \eps \circ f$.
\end{lemma}
\begin{proof}
If $f: C \to C$ is right $C$-colinear and left $D$-colinear,
then for all $c \in C$, $f(c\1) \o c\2 = \cop(f(c))$ and $g(c\1) \o f(c\2) = (g \o C)\cop(f(c))$,
whence $$ g(c\1) \eps(f(c\2)) = g(f(c)) = \eps(f(c\1)) g(c\2). $$
It follows that $g^*(d^*) (\eps \circ f) = (\eps \circ f) g^*(d^*)$ for
all $d^* \in D^*$ w.r.t.\ the convolution product. 

The inverse mapping $  V_{C^*}(D^*) \to \End {}^DC^C$ is given by
$$ c^* \longmapsto (c \mapsto c^*(c\1)c\2). $$
We obtain $g(c\1) c^*(c\2) = c^*(c\1) g(c\2) $ (for each $c \in C, c^* \in V_{C^*}(D^*)$)
since there is equality when any $d^* \in D^*$ is applied (and any $K$-dual of coalgebras
we consider ``separates points'').  It follows that the mapping $c \mapsto c^*(c\1)c\2$
is left $D$-colinear. 
Of course, by right $C$-colinearity of $f$ we have $\eps(f(c\1))c\2 = f(c)$
for each $c \in C$. 
It is easy to see
that we have defined an anti-isomorphism of algebras.  
\end{proof}

Now suppose $M$ is a $D$-$C$-bicomodule where $g: C \to D$ continues to be a coalgebra homomorphism.
It is well-known that $M$ then also has $C^*$-$D^*$-bimodule structure via  convolution actions.  From this we define the $D^*$-$C^*$-bimodule structure on the $K$-dual $M^*$ by $(d^* \cdot m^* \cdot c^*)(m) = d^*(m\-1) m^*(m\0) c^*(m\1)$.
Similar to the lemma we prove:

\begin{prop}
If $g : C \to D$ is a coalgebra homomorphism and $M$ is a $D$-$C$-bicomodule,
then ${\rm Hom}^{D-C}(M,C) \cong (M^*)^{D^*}.$
\end{prop}

\begin{cor}
Under the conditions above, we have
${\rm Hom}^{D-C}(C \b_D C, C) \stackrel{\cong}{\longrightarrow} (C \b_D C)^{* \, D^*}$
via $ f \mapsto \eps \circ f$.
\end{cor}
\begin{proof}
We note that the inverse mapping is given by
\begin{equation}
\label{eq: map2}
\eta \longmapsto (c \o c' \mapsto \eta(c \o {c'}\1) {c'}\2) \qed
\end{equation} 
\renewcommand{\qed}{}\end{proof}

\subsection{Left coD2 quasibases.}  We are now in a position to re-write eq.~(\ref{eq: basic}) identifying
$f_i \in {\rm Hom}^{D-C}(C, C\b_DC)$ with
$\alpha_i \in \End {}^DC^D$ using the mapping~(\ref{eq: map1}), and identifying $g_i \in {\rm Hom}^{D-C}(C \b_D C, C)$ with $\eta_i \in (C \b_D C)^{* \, D^*}$
using mapping~(\ref{eq: map2}).  We obtain for each
$c \o c' \in C \b_D C$:
\begin{equation}
\label{eq: lcd2qb}
 c \o c' = \sum_{i=1}^N \eta_i(c \o {c'}\1) \alpha_i({c'}\2) \o {c'}\3.
\end{equation}
The equation is analogous to the eq.~(\ref{eq: lqbd2}); for that reason we call $\eta_i \in (C \b_D C)^{*\, D^*}$
and $\alpha_i \in \End {}^DC^D$ \textit{left coD2
quasibases} for the coalgebra homomorphism $g: C \to D$. The quasibase equation above has the equivalent form, 
\begin{equation}
\label{eq: diagrammatic}
C \, \b_D \, C = \sum_{i=1}^N (\eta_i \o \alpha_i \o C) \circ (C \o \cop^2)
\end{equation}
 
Given $\beta \in \End {}^DC^D$, we note that
for $c \in C$ we have $\beta(c\1) \o c\2 \in C \b_D C$,
and substituting this into eq.~(\ref{eq: lcd2qb})
and applying $C \o \eps$ to this yields
\begin{equation}
\label{eq: lemma}
\beta(c) = \sum_{i=1}^N \eta_i(\beta(c\1) \o c\2)\alpha_i(c\3),
\end{equation}
in other words, $\beta = \sum_i (\eta_i \o \alpha_i)\circ (\beta \o C \o C)\circ \cop^2$. 

\subsection{Right bialgebroid structure on $\End {}^DC^D$ over $C^{* \, D^*}$.} We proceed to show
that a coD2 coalgebra homomorphism $g: C \to D$ has
bialgebroid structure on $\End {}^DC^D$.  This and its  noncommutative base algebra will be denoted by  
\begin{equation}
 E := \End {}^DC^D,  \ \ \ R := C^{* \, D^*}.
\end{equation}
There are immediately two commuting mappings, a
homomorphism $s: R \to E$ and an anti-homomorphism
$t: R \to E$, given by ($r \in R, c \in C$)
\begin{equation}
\label{eq: source}
s(r)(c) = c\1 r(c\2),
\end{equation}
a source map, and a target map, 
\begin{equation}
\label{eq: target}
t(r)(c) = r(c\1) c\2.
\end{equation}
We note that $s(r) t(r') = t(r') s(r)$ w.r.t.\
composition in $E$ since both applied to $c \in C$
yield $r'(c\1) r(c\3) c\2$ by coassociativity of $\cop$.
We note that the $R$-$R$-bimodule structure induced
by $s$ and $t$ from the right, suggestively denoted
by $E_{s,t}$ is given  by the straightforward 
\begin{equation}
\label{eq: bimodule}
(r \cdot \alpha \cdot r' )(c) = r(c\1) \alpha(c\2) r'(c\3)
\end{equation}  
for $r,r' \in R, \alpha \in E, c \in C$.  

At this point we may profitably note that eq.~(\ref{eq: lemma})
shows that ${}_RE$ is f.g.\ projective, since the $N$ mappings 
$\beta \mapsto \eta_i \circ (\beta \o C) \circ \cop$
are easily seen to be in $\Hom ({}_RE, {}_RR)$.  

The counit map $\eps_E: E \to R$ is given by
\begin{equation}
\eps_E(\alpha) = \eps \circ \alpha
\end{equation}
where $\eps: C \to K$ denotes the counit on $C$.
Since  $\alpha$ in $E$ is right and left $D$-colinear,
it follows that $g(c\1) (\eps \circ \alpha)(c\2) = 
(\eps \circ \alpha)(c\1) g(c\2)$ for each $c \in C$ since both equal
$g(\alpha(c))$.  Then of course $\eps \circ \alpha
\in C^{* \, D^*} = R$.  Note that $\eps_E(\id_C) = \eps = 1_R$.  Note the $\eps_E$ is left and right $R$-linear:
$$ \eps_E(r \cdot \alpha \cdot r')(c) = r(c\1) \eps(\alpha(c\2)) r'(c\3) = (r \eps_E(\alpha) r')(c) $$
w.r.t.\ the convolution product in $R$.  Also
note that $(\alpha, \beta \in E, c \in C$)
\begin{eqnarray*}
\eps_E(s(\eps_E(\alpha)) \circ \beta)(c) & = & \eps \circ s(\eps_E(\alpha))(\beta(c)) \\
& = & \eps({\beta(c)}\1) \eps(\alpha({\beta(c)}\2)) \\
& = & \eps_E(\alpha \circ \beta)(c) 
\end{eqnarray*}
Thus, $\eps_E(\alpha \circ \beta) = \eps_E(s(\eps_E(\alpha)) \circ \beta)$ and similarly
we compute $\eps_E(\alpha \circ \beta) = \eps_E(t(\eps_E(\alpha)) \circ \beta)$.
  
Dualizing the definition of generalized
Lu coproduct in \cite[eq.\ (74)]{KS}, we would obtain a coproduct on $E$
given by $\cop_E(\alpha) = \cop \circ \alpha$ in ${\rm Hom}^{D-D}(C, C \b_D C)$.  In order to 
make sense of this, we need:

\begin{prop}
\label{prop-main}
If $g: C \to D$ is left coD2, then there is an isomorphism, 
\begin{equation}
\label{eq: map0}
E \o_R E \stackrel{\cong}{\longrightarrow} {\rm Hom}^{D-D} (C, C \b_D C), \ \ \ \alpha \o \beta \longmapsto ( c \mapsto \alpha(c\1) \o \beta(c\2))
\end{equation}
\end{prop}
\begin{proof}
It is clear from $D$-colinearity of $\alpha$ and $\beta$
as well as eq.~(\ref{eq: bimodule}) 
that the mapping $a \o_R \beta \mapsto (\alpha \o \beta)\circ \underline{\cop}$ is well-defined in all
respects.  

An inverse mapping ${\rm Hom}^{D-D}(C, C \b_D C) \to E \o_R E$ is given in terms of the left coD2 quasibases $\eta_i \in (C \b_D C)^{* \, D^*}$ and $\alpha_i \in E$ in eq.~(\ref{eq: lcd2qb}): 
\begin{equation}
\label{eq: identify}
G \longmapsto \sum_{i=1}^N (C \o \eta_i) \circ (G \o C) \circ \underline{\cop} \o_R \alpha_i,
\end{equation}
where $G \in {\rm Hom}^{D-D}(C, C \b_D C)$.  
It is easy to see that $(C \o \eta_i) \circ (G \o C) \circ \underline{\cop}$ is left $D$-colinear,
but it is right $D$-colinear as well, since ($c \in C$)
$$(C \o g) \cop (C \o \eta_i)(G(c\1) \o c\2) = (C \o \eta_i)(G(c\1) \o c\2) \o g(c\3) $$
follows from $G$ taking values in $C \b_D C$ and $\eta_i$ satisfying $g(c\1) \eta_i(c\2 \o c')
= \eta_i(c \o {c'}\1) g({c'}\2)$ for $c \o c' \in C \b_D C$.   
By eqs.~(\ref{eq: lemma}) and~(\ref{eq: bimodule}) this mapping sends
$$ (\alpha \o \beta) \circ \underline{\cop} \mapsto 
\sum_{i=1}^N \alpha \o_R (\eta_i \circ (\beta \o C) \circ \underline{\cop})\cdot \alpha_i = \alpha \o \beta,$$
where $\eta_i \circ (\beta \o C) \circ \underline{\cop} \in R$.  

Conversely, note that the image of $G \in {\rm Hom}^{D-D}(C, C \b_D C)$ in $E \o_R E$ is sent into
$$ \sum_{i=1}^N (C \o \eta_i \o C) \circ (G \o C \o \alpha_i) \circ \underline{\cop}^2 = G,$$
since $\id_{C \b_D C} = \sum_{i=1}^N (\eta_i \o \alpha_i \o C) \circ (C \o \underline{\cop}^2)$. 
\end{proof}

Note that $C$ has left $R$-module structure that commutes with the right coaction $\rho^R:
C \to C \o D$.  The left module action is naturally $r \cdot c = c\1 r(c\2)$ where $r \in
R = C^{*\, D^*}$. Then note that 
$$ \rho^R(r \cdot c) = c\1 \o g(c\2) r(c\3) =  c\1 r(c\2) \o g(c\3) = r \cdot \rho^R(c).$$
Similarly, ${}^DC_R$ has associative action and coaction.  
It follows that ${\rm Hom}^{D-D}(C, C \b_D C)$ has $R$-$R$-bimodule structure induced by its contravariant argument:  
\begin{equation}
\label{eq: are}
(r \cdot G \cdot r')(c) = r(c\1) G(c\2) r'(c\3) \ \ \ (r \in R, c \in C, G \in {\rm Hom}^{D-D}(C, C \b_D C) )
\end{equation}
Also $E \o_R E$ has a natural $R$-$R$-bimodule structure derived from ${}_RE$ and $E_R$
in the first and second tensorands.  The isomorphism in the proposition is clearly
left and right $R$-linear:
\begin{cor}
The isomorphism $E \o_R E \stackrel{\cong}{\longrightarrow} {\rm Hom}^{D-D} (C, C \b_D C)$,
given by \newline
$\alpha \o_R \beta \mapsto (\alpha \o \beta) \circ \underline{\cop}$,
is an $R$-$R$-bimodule isomorphism.
\end{cor}

By identifying $E \o_R E$ with ${\rm Hom}^{D-D}(C, C \b_D C)$, we define a comultiplication for $E$ by 
\begin{equation}
\cop_E(\alpha) = \underline{\cop} \circ \alpha,
\end{equation}
clearly a right and left $D$-colinear homomorphism from $C$ into $C \b_D C$
for each $\alpha \in E = \End {}^DC^D$.  We note that $\cop_E$ is $R$-$R$-linear:
$$ \cop_E(r \cdot \alpha \cdot r')(c) = r(c\1) \underline{\cop}(\alpha(c\2)) r'(c\3) =
(r\cdot \cop_E(\alpha)\cdot r')(c), $$
by eq.~(\ref{eq: are}).  

Note that $\cop_E(1_E) = 1_E \o_R 1_E$, since
$\cop_E(\id_C) = \underline{\cop}$, which corresponds
to $\id_C \o \id_C$ under the identification mapping~(\ref{eq: map0}). 

Next we show that $(\eps_E \o_R E) \circ \cop_E = E$ (where $E$ denotes $\id_E$).  Note that via the identification
mapping~(\ref{eq: identify}) the formal definition of the comultiplication $\cop_E:
E \to E \o_R E$ is ($\beta \in E$)
\begin{equation}
\label{eq: formal}
\cop_E(\beta) = \sum_{i=1}^N (C \o \eta_i) \circ ((\underline{\cop} \circ \beta) \o C) \circ \underline{\cop} \o \alpha_i.
\end{equation}
Whence identifying $R \o_R E \cong E$ canonically and letting $c \in C$: 
$$ (\eps_E \o_R E) \circ \cop_E (\beta)(c) = \sum_{i=1}^N \eps({\beta(c\1)}\1)\eta_i({\beta(c\1)}\2 \o c\2)
\alpha_i(c\3) $$
$$ = \sum_i \eta_i(\beta(c\1) \o c\2) \alpha_i(c\3) = \beta(c),$$
by eq.~(\ref{eq: lemma}).   

Next note that $(E \o_R \eps_E) \circ \cop_E = E$ as follows.  Apply $\cop \o C$
to the quasibases eq.~(\ref{eq: lcd2qb}) for $c \o c' \in C \b_D C$ to obtain
$$ c\1 \o c\2 \o c' = \sum_{i=1}^N c\1 \eta_i(c\2 \o {c'}\1) \o \alpha_i({c'}\2) \o {c'}\3, $$
to which we apply $C \o \eps \o \eps$, obtaining
\begin{equation}
\label{eq: star}
c \eps(c') = \sum_i c\1 \eta_i(c\2 \o {c'}\1) \eps(\alpha_i({c'}\2).
\end{equation}
We apply this to $\beta(c\1) \o c\2 \in C \b_D C$ to see that
\begin{eqnarray*}
((E \o \eps_E) \circ \cop_E)(\beta)(c) & = & \sum_i {\beta(c\1)}\1 \eta_i({\beta(c\1)}\2 \o c\2)
\eps(\alpha_i(c\3)) \\
& = & \beta(c\1) \eps(c\2) = \beta(c). 
\end{eqnarray*}

Denote the values $\cop_E(\beta) = \beta\1 \o_R \beta\2 \in E \o_R E$ using Sweedler's notation.
In ${\rm Hom}^{D-D} (C, C\b_D C)$ this is identified with the mapping
$c \mapsto \beta\1(c\1) \o \beta\2(c\2)$.  At the same time, we define
$\cop_E(\beta) = \underline{\cop} \circ \beta$ in this same hom-group, thus for
every $\beta \in E$ and $c \in C$, 
\begin{equation}
\label{eq: distribution}
{\beta(c)}\1 \o {\beta(c)}\2 = \beta\1(c\1) \o \beta\2(c\2),
\end{equation}
as an equality in $C \o_K C$.

Now note that for each $r \in R$ and $\alpha \in E$, we have
\begin{eqnarray*}
(s(r)\circ \alpha\1 )(c\1) \o \alpha\2 (c\2) & = & s(r)({\alpha(c)}\1) \o {\alpha(c)}\2 \\
& = & {\alpha(c)}\1 r({\alpha(c)}\2) \o \alpha(c)\3 \\
& = & {\alpha(c)}\1 \o t(r)({\alpha(c)}\2) \\
& = & \alpha\1(c\1) \o (t(r) \circ \alpha\2)(c\2), 
\end{eqnarray*} 
whence $\Im \cop_E \subseteq (E \o_R E)^R$ where the $R$-$R$-bimodule structure
on ${}_{\cdot}E \o_R E_{\cdot}$ uses ${}_{s, t}E$ instead.  

As is well-known in the theory of bialgebroids, the last condition on $\Im \cop_E$
shows that there is a well-defined tensor algebra multiplication on $\Im \cop_E$.  
We claim that
\begin{equation}
\cop_E(\alpha \circ \beta) = \alpha\1 \circ \beta\1 \o_R \alpha\2 \circ \beta\2
\end{equation}
Using the identification $E \o_R E \cong {\rm Hom}^{D-D}(C, C \b_D C)$ in proposition~\ref{prop-main} again,
for $c \in C$,
\begin{eqnarray*}
\alpha\1(\beta\1(c\1)) \o \alpha\2(\beta\2(c\2)) & = & \alpha\1({\beta(c)}\1) \o \alpha\2({\beta(c)}\2) \\
{\alpha(\beta(c))}\1 \o {\alpha(\beta(c))}\2 & = &  (\alpha \circ \beta)\1 (c\1) \o (\alpha \circ \beta)\2 (c\2) ,
\end{eqnarray*}
hence $\cop_E(\alpha \circ \beta) = \cop_E(\alpha) \cop_E(\beta)$ in this
special submodule of $E \o_R E$ where tensor algebra multiplication is valid.  

Finally we note that $\cop_E$ is coassociative.
Heuristically this depends on the coassociativity of 
the coproduct $\cop$ on $C$ since $\cop_E(\alpha) =
\underline{\cop} \circ \alpha$ where $\underline{\cop}$
is just $\cop$ with a restriction of its codomain.
To prove coassociativity we need a lemma
generalizing proposition~\ref{prop-main}:

\begin{lemma}
Suppose $M$ is a $C$-$C$-bicomodule and $g : C\to D$
a coalgebra homomorphism of codepth two.  
Then there is a natural $R$-$R$-bimodule isomorphism,
$${\rm Hom}^{D-D}(C, M) \o_R E \cong {\rm Hom}^{D-D}(C, M \b_D C), $$
via $f \o_R \alpha \longmapsto (c \mapsto
 f(c\1)\o \alpha(c\2))$.
\end{lemma}
\begin{proof} The proof is almost identical to the proof
of proposition~\ref{prop-main}, and is therefore omitted.  
\end{proof}
 
Applying the lemma  with $M = C \b_D C$,  proposition~\ref{prop-main} and using coassociativity of cotensor product \cite[11.6]{BW}, 
we note that
\begin{equation}
\label{eq: three}
E \o_R E \o_R E  \cong  {\rm Hom}^{D-D} (C, C \b_D C \b_D C)
\end{equation}
 via $$\alpha \o_R \beta \o_R \gamma \longmapsto ( c \mapsto \alpha(c\1) \o \beta(c\2) \o \gamma(c\3) ) .$$
It follows from eq.~(\ref{eq: distribution}) that  
$(\cop_E \o E)\cop_E(\alpha)$ is identified with
$${\alpha\1 (c\1)}\1 \o {\alpha\1(c\1)}\2 \o \alpha\2(c\2) =
\alpha_{(1,1)}(c\1) \o \alpha_{(1,2)}(c\2) \o \alpha\2(c\3), $$
but the LHS equals
$$ \alpha\1(c\1) \o {\alpha\2(c\2)}\1 \o \alpha\2(c\2)\2 = \alpha_1(c\1) \o \alpha_{(2,1)}(c\2)
\o \alpha_{(2,2)}(c\3)  $$
which is identified with $= (E \o_R \cop_E)\cop_E(\alpha)$. Via the isomorphism~(\ref{eq: three}) we obtain $(\cop_E \o_R E) \circ \cop_E = (E \o_R \cop_E) \circ \cop_E$.

We have proven:
\begin{theorem}
If $g: C \to D$ is a left codepth two homomorphism of coalgebras,
then $E = \End {}^DC^D$ is a right bialgebroid over $R = C^{* \, D^*}$,
the centralizer subalgebra in $C^*$ induced by the dual algebra homomorphism $g^*: D^* \to C^*$.
Moreover, the left $R$-module $E$ is finitely generated projective.  
The structure of the  $R$-bialgebroid  and endomorphism algebra $E$ is given by ($r,r' \in R, c \in C, \alpha \in E$)
\begin{eqnarray}
s(r)(c) & = & c\1 r(c\2) \\
t(r)(c) & = & r(c\1) c\2 \\
(r \cdot \alpha \cdot r')(c) & = & r(c\1) \alpha(c\2) r'(c\3) \\
\cop_E(\alpha) & = & \sum_{i=1}^N (C \o \eta_i) \circ ((\underline{\cop} \circ \alpha) \o C) \circ \underline{\cop} \o_R \alpha_i \\
\eps_E(\alpha) & = & \eps \circ \alpha 
\end{eqnarray}
with respect to the codepth two quasibases $\alpha_i \in E$ and $\eta_i \in (C \b_D C)^{* \, D^*}$
satisfying for all $c \o c' \in C \b_D C$, 
\begin{equation} 
 c \o c' = \sum_{i=1}^N \eta_i(c \o {c'}\1) \alpha_i({c'}\2) \o {c'}\3.
\end{equation}
\end{theorem}
\begin{proof}
For the convenience of the reader, we gather the axioms of a right $R$-bialgebroid $E$,
verified in the subsection preceding the theorem.  Given $K$-algebras $R$ and $E$,
a bialgebroid $E$ over $R$ has ``source'' algebra homomorphism $s: R \to E$ and ``target''
algebra anti-homomorphism $t: R \to E$ such that commutativity $s(r) t(r') = t(r')s(r)$ 
holds within $E$ for all $r,r' \in R$.  We then refer to an $R$-$R$-bimodule structure
on $E$ induced by $r \cdot e \cdot r' = e s(r') t(r)$ for all $r,r' \in R, e \in E$,
and the natural $R$-$R$-bimodule structures on $R$ and $E \o_R E$.  
The
axioms are then:
\begin{enumerate}
\item There is an $R$-coring $(E,R,\cop_E, \eps_E)$: i.e., ``comultiplication'' $\cop_E: E \to E\o_R E$ and ``counit'' $\eps_E: E \to R$ are right
and left $R$-linear, $\cop_E$ is coassociative, and $\eps_E$ satisfies counitality axioms;
\item The comultiplication and counit are unit-preserving:
$\cop_E(1_E) = 1_E \o_R 1_E$ and $\eps_E(1_E) = 1_R$;
\item the comultiplication takes its values in the submodule of finite sums,
$$ E \times_R E := \{ \sum_i x_i \o_R y_i \in E \o_R E \, | \, \forall \, r \in R, 
\sum_i s(r) x_i \o_R y_i = \sum_i x_i \o_R t(r)y_i \} $$
where tensor algebra multiplication is well-defined;
\item for all $e,e' \in E$, we have $\cop_E(ee') = \cop_E(e) \cop_E(e') $;
\item the unital tensor category axiom for $R$ as an $E$-module: $$\eps_E(ee') = \eps_E(s(\eps_E(e)) e') = \eps_E(t(\eps_E(e))e')$$ for all $e,e' \in E$. \qed
\end{enumerate}
\renewcommand{\qed}{}\end{proof}

\subsection{Discussion.}
A theory of right coD2 coalgebra homomorphisms has a 
similar development. By consulting \cite{BW}, one might under suitable hypotheses develop
a similar theory of codepth two for a homomorphism of $R$-corings.
The (``co-Sweedler'') $C$-ring 
$C \b_D C$ for any coalgebra homomorphism $g: C \to D$ is defined by $\mu: C \b_D C \b_D C \to C \b_D C$, $\mu(c \o c' \o c'') = c \eps(c') \o c''$
with unit $\eta = \underline{\cop}: C \to C \b_D C$.
Something related to this $C$-ring might play a role
in a more complete theory of codepth two.   

If the coalgebra homomorphism $g: C \to D$ is coD2
and its dual $g^*: D^* \to C^*$ is D2,
 then there
is an anti-monomorphism of $R$-bialgebroids
$E \to S := \End {}_{D^*}C^*_{D^*}$ given by
$\alpha \mapsto \hat{\alpha}$, where $\hat{\alpha}(c^*)
= c^* \circ \alpha$.  If $C$ and $D$ are finite dimensional, this is an anti-isomorphism of a left and
right bialgebroid over $R$.  It would be interesting to know something of  the precise
relationship between D2 and coD2.  For example,
 the quotient homomorphism $H \to \overline{H}$ should be coD2 if $K \into H$ is a normal
Hopf subalgebra with $\overline{H} = H / HK^+$.  

Given a coD2 coalgebra homomorphism $g: C \to D$
and its constructions $E$ and $R$ defined above,
the coalgebra $C$ is a left $E$-module coalgebra
(i.e. a coalgebra in the tensor category of left $E$-modules) due to eq.~(\ref{eq: distribution})
and one other, counital axiom. If ${}^DC$ is
co-balanced, it should be  that
$\ker g = \{ \alpha(c) - c\1 \eps(\alpha(c\2)) |
c \in C, \alpha \in E \}$ and we obtain the image
of $g$ in $D$ in a type of coGalois theory.     

It would be interesting to pursue the rest of the
structure of duality, e.g. endomorphism algebras
and smash products, realization of the $R$-dual
left bialgebroid of $E$, and its relationship
to \cite{CQR}.

%%%%%%%%%%%%%%%%%%%%%%%%%%%%%%%%%%%%%%%%%%%%%%%%%%%%%%%%%%%%%%%%%%%%%%%%%%%%%%%%%%%%%%%%%%%%%%%%

\end{document}